\documentclass[12pt,oneside,reqno]{amsart}
\usepackage{graphicx}
\usepackage{caption,subcaption}
\usepackage{dcolumn}
\usepackage{amsmath,amsthm,amssymb}
\usepackage{todonotes}
\usepackage{url}
\usepackage{placeins}

\usepackage{array,booktabs}

\renewcommand{\Re}{\operatorname{Re}}

\newtheorem*{theorem*}{Theorem}
\newtheorem*{conj*}{Conjecture}
\theoremstyle{remark}
\newtheorem*{rem*}{Remark}

\begin{document}

\title[Threshold of blowup for equivariant wave maps]{Threshold for blowup for equivariant wave maps in higher dimensions}

\author{Pawe\l {} Biernat}
\address{Mathematisches Institut, Universit\"at Bonn, Germany }
\email{pawel.biernat@gmail.com}
\author{Piotr Bizo\'n}
\address{Institute of Physics, Jagiellonian
University, Krak\'ow, Poland\\
and
Max Planck Institute for Gravitational Physics (Albert Einstein Institute),
Golm, Germany}
\email{bizon@th.if.uj.edu.pl}
\author{Maciej Maliborski}
\address{
Max Planck Institute for Gravitational Physics (Albert Einstein Institute),
Golm, Germany}
\email{maciej.maliborski@aei.mpg.de}

\date{\today}%
\begin{abstract}
 We consider equivariant wave maps from $\mathbb{R}^{d+1}$ to $\mathbb{S}^d$  in supercritical dimensions $3\leq d\leq 6$. Using mixed numerical and analytic methods, we show that the threshold of blowup is given by the codimension-one stable manifold of a self-similar solution with one instability.
\end{abstract}
\maketitle

\section{Introduction}
This paper is concerned with the wave map equation
\begin{equation}\label{wmeq}
  \phi_{tt}-\Delta \phi + \left(|\phi_t|^2-|\nabla \phi|^2 \right)\, \phi=0\,,
\end{equation}
where $(t,x)\in  \mathbb{R} \times \mathbb{R}^d$ and $\phi(t,x)\in \mathbb{S}^d \hookrightarrow \mathbb{R}^{d+1}$.
We assume that $d\geq 3$ and restrict attention to equivariant maps of the form (where $r=|x|$)
\begin{equation}\label{ansatz}
  \phi(t,x)=\left(\frac{x}{r}\,\sin{u(t,r)},\cos{u(t,r)}\right)\,.
\end{equation}
For this ansatz Eq.\eqref{wmeq} reduces to the scalar semilinear wave equation
\begin{equation}\label{eqv}
  u_{tt} = u_{rr} +\frac{d-1}{r}\, u_r -\frac{d-1}{2 r^2}\, \sin(2u)\,.
\end{equation}
The goal is to understand global dynamics for smooth initial data $(u,u_t)\vert_{t=0}$.
In \cite{bb} we found  the explicit stable self-similar solution \eqref{f0}  and gave numerical and analytic evidence that for all $d\geq 3$ this solution determines the universal asymptotics of generic blowup for large initial data (for $d=3$ this was established earlier in \cite{bct} and \cite{d1}). On the other hand, it is well known that solutions starting from small initial data remain globally regular in time. The dichotomy between global regularity and blowup raises a natural question about the nature of a borderline between these two generic asymptotic behaviors. This question was first studied for $d=3$ in \cite{bct} which gave evidence that the threshold for blowup is determined by the codimension-one stable manifold of a self-similar solution with one instability
(whose existence was established in \cite{b1}).
In this paper we extend this analysis to dimensions $4\leq d\leq 6$. In higher dimensions  the threshold dynamics  is qualitatively different and will be described elsewhere.

The outline of the paper is as follows. In section~2 we provide classification of self-similar solutions of Eq.\eqref{eqv}  and in section~3 we analyze their linear stability. Finally, in section~4 we present  results of numerical computations of dynamics at the threshold for blowup.

\section{Self-similar solutions}
By definition, self-similar solutions are invariant under the scaling,  $u(t/L,r/L)=u(t,r)$, hence they have the form
\begin{equation}\label{css-ansatz}
  u(t,r)=f(y)\,,\qquad y=\frac{r}{T-t}\,,
\end{equation}
where a positive constant $T$, clearly allowed by the time translation symmetry, is introduced for later convenience.
Inserting this ansatz into Eq.\eqref{eqv} we obtain the ordinary differential equation
\begin{equation}\label{css}
(1-y^2) f'' +\left(\frac{d-1}{y}-2y\right) f'
-\frac{d-1}{2y^2}\,\sin(2f)=0\,.
\end{equation}
We require  solutions to be smooth on closed interval $0\leq
y\leq 1$, which corresponds to the solid past light
cone of the  point $(t=T,r=0)$. For such solutions
\begin{equation}\label{css-blow}
  \partial_r  f\left(\frac{r}{T-t}\right )\Big\vert_{r=0} = \frac{f'(0)}{T-t} \rightarrow \infty \quad\mbox{as}\,\,\, t\nearrow T\,,
\end{equation}
hence each self-similar solution $f(y)\in C^{\infty}[0,1]$ is an example of a singularity developing in finite time from smooth initial data.

It follows from \eqref{css} that
for local smooth solutions near the origin
\begin{equation}\label{y=0}
f(y)= c y+\mathcal{O}(y^3)\,,
 \end{equation}
 where the parameter $c$ determines the solution uniquely. It is not difficult to  show  that these local solutions can be extended smoothly to the whole interval $0\leq y<1$ \cite{b1} but, in general, they are not smooth at $y=1$. The classification of self-similar solutions amounts to finding all (isolated) values of $c$ for which the solution \eqref{y=0} is smooth at $y=1$. One such value is $c_0=\frac{2}{\sqrt{d-2}}$ for which the solution is known is closed form \cite{bb} (for $d=3$ this solution was known earlier  \cite{s, ts})
 \begin{equation}\label{f0}
  f_0(y)=2 \arctan\left(\dfrac{y}{\sqrt{d-2}}\right)\,.
  \end{equation}
  We conjecture that $f_0$ is the only self-similar solution for $d\geq 7$.
  To find other solutions in dimensions $3\leq d\leq 6$, let us note that if $f(y)$ is smooth at $y=1$, then \cite{cst}
  \begin{subequations}
  \begin{align}
    (d-3) f'(1) &-\frac{d-1}{2} \sin(2f(1))=0, \\
    (d-5) f''(1) &+ \left(d-7-(d-1)\cos(2f(1))\right) f'(1)=0,
  \end{align}
  \end{subequations}
 where (9a) follows directly from \eqref{css} and (9b) follows from multiplying \eqref{css} by $y^2$ and differentiating.  As a consequence, smooth solutions have the following Taylor series expansions at $y=1$:
 \begin{itemize}
 \addtolength{\itemsep}{0.1\baselineskip}
 \item For $d=3$ we have
 \begin{equation}\label{d3}
   \hspace{-2.6cm}f(y)=\frac{\pi}{2}-f'(1) (1-y) +....
 \end{equation}
 where $f'(1)$ is the only free parameter.
\item  For $d=5$ there are two possibilities. Either
 \begin{equation}\label{case1}
   \hspace{-2.6cm} f(y)=\frac{\pi}{2} + \frac{1}{2} f''(1) (1-y)^2  +...
 \end{equation}
 or
 \begin{equation}\label{case2}
   f(1)=\frac{\pi}{3}-\frac{\sqrt{3}}{2}\, (1-y)+ \frac{1}{2} f''(1) (1-y)^2+....
 \end{equation}
 In both cases $f''(1)$ is the only free parameter.
 \item For $d=4,6$ we have
  \begin{equation}\label{even}
   f(y)=f(1) - \frac{d-1}{2(d-3)} \,\sin(2f(1)) (1-y)  +...
 \end{equation}
 where $f(1)\neq \pi/2$ is the only free parameter.
 \end{itemize}
 \begin{theorem*}
For each $d\in\{3,4,5,6\}$ there is an infinite sequence $(c_n)_{n\in \mathbb{N}}$ such that the  solution
\eqref{y=0} is smooth at $y=1$ and behaves as \eqref{d3} for $d=3$, as \eqref{case1} for $d=5$, and as \eqref{even} for $d=4,6$.  We denote the corresponding solutions by $f_n(y)$.
\end{theorem*}
For $d=3,5$ this theorem was proven in \cite{b1} using
 a shooting argument. Key to this argument is  the change of an independent variable
which brings equation \eqref{css} into an asymptotically autonomous form that can be analyzed by dynamical system methods. For $d=4,6$ the proof requires a minor modification that we leave to the interested reader as an exercise. For $d=5$  an alternative  variational proof of existence of $f_1(y)$ was given  in \cite{cst}.

The index $n$ on $f_n$ is the nodal index. From the construction of self-similar solutions given in  \cite{b1} it follows that
\begin{equation}\label{N}
n=\begin{cases}
N(n,d)\quad & \mbox{for} \,\, d=3,4,\\
N(n,d)+1\quad & \mbox{for} \,\, d=5,6\,,
\end{cases}
\end{equation}
where $N(n,d)$ be the number of zeros of $f_n'(y)$ on $(0,1)$.

\section{Linear stability analysis}
 As the first step towards understanding the role of the  self-similar solutions $f_n(\frac{r}{T-t})$ in dynamics we need to analyze their linear stability. To this end
 it is convenient to define a new time coordinate
$s=-\ln(T-t)$ and rewrite Eq.\eqref{eqv} in terms of
$U(s,y)=u(t,r)$
\begin{equation} \label{eq-ys}
U_{ss} + U_{s} + 2 y\: U_{sy}
=(1-y^2) U_{yy} +\left(\frac{d-1}{y}-2y\right) U_{y}
-\frac{d-1}{2y^2}\,\sin(2 U)\,.
\end{equation}
In these variables the problem of finite time self-similar blowup in converted
into the problem of asymptotic convergence for $s \rightarrow
\infty$ towards a stationary solution $f(y)$.
 Following the standard procedure we seek
solutions of Eq.\eqref{eq-ys} in the form
$U(s,y)=f_n(y)+ e^{\lambda s} v(y)$. Dropping nonlinear terms
 we get the quadratic eigenvalue problem  on the interval $0\leq
y\leq 1$
\begin{equation}\label{eigen}
(1-y^2) v''+\left(\frac{d-1}{y}-2(\lambda+1)y\right) v' -\lambda(\lambda+1) v - \frac{d-1}{y^2}\,\cos(2 f_n)\, v=0.
\end{equation}
 By assumption,  the solution $U(s,y)$ is smooth for $s<\infty$, hence
 we demand that  $v\in \mathcal{C}^{\infty}[0,1]$. This condition leads to the quantization of the eigenvalues. The Frobenius indices  for $v(y)$ are $\{1,1-d\}$ at $y=0$ and $\{0,\frac{d-1}{2}-\lambda\}$ at $y=1$,  where the first index in each pair gives the smooth solution.

The linear stability analysis of the solution $f_0$ in \cite{bb, cdg} took advantage of the fact that in this case Eq.\eqref{eigen} can be expressed as the Heun equation. Unfortunately, other self-similar solutions are not known in closed form so this approach is not possible.

Let us denote the eigenvalues and eigenfunctions for $f_n$ by $\lambda_k^{(n)}$ and  $v_k^{(n)}(y)$.
    We conjecture that for each $n$  the eigenvalues are real and can be ordered as follows
    \begin{equation}\label{spectrum}
   \underbrace{\dots < \lambda_{-2}^{(n)}<\lambda_{-1}^{(n)}}_\text{$\infty$ many stable modes}<0<\underbrace{\lambda_{0}^{(n)}=1}_\text{gauge mode}<\underbrace{\lambda_{1}^{(n)}<\dots<\lambda_{n}^{(n)}}_\text{$n$ unstable modes}\,.
   \end{equation}
The eigenvalue $\lambda_0^{(n)}$ corresponds to the gauge mode $v_0^{(n)}(y)=y\, f'_n(y)$ generated by the shift of the blowup time $T$. In the following we corroborate \eqref{spectrum} by  analytic and numerical arguments.

It is instructive to consider a self-adjoint eigenvalue problem associated with \eqref{eigen}.
Let
\begin{equation}\label{w}
w(y)=(1-y^2)^{\lambda/2}\, y^{\frac{d-1}{2}}\, v(y)\,.
\end{equation}
 Then, Eq.\eqref{eigen}  becomes
\begin{equation}\label{An}
  A_n w = \mu w, \qquad A_n =-(1-y^2)^{\frac{d+1}{2}} \partial_y \left((1-y^2)^{\frac{3-d}{2}} \partial_y\right) + V_n(y)\,,
 \end{equation}
where
\begin{equation}\label{mu}
\mu=\lambda(d-1-\lambda)
\end{equation}
and
$$
V_n(y)=\frac{d-1}{y^2} (1-y^2)  \left(\cos(2f_n(y)) +\frac{(d-3)}{4}(1+y^2)\right) \,.
$$
The operator $A_n$ is essentially self-adjoint on the Hilbert space $X=L^2\left([0,1],(1-y^2)^{-\frac{d+1}{2}}\,dy\right)$. We shall refer to the eigenvalues of $A_n$ as $\mu$-eigenvalues, not to be confused with the eigenvalues $\lambda$ of our problem.

The Frobenius indices for $w(y)$ are $\frac{1}{2}(1\pm d)$ at $y=0$ and $\frac{1}{4}(d-1\pm\sqrt{(d-1)^2-4\mu})$ at $y=1$, where  the $+$ indices  give solutions belonging to $X$ (for $\mu<\frac{1}{4}(d-1)$, the bottom of the continuous spectrum). For $\lambda>\frac{d-1}{2}$ we have
$$
\frac{1}{4} \left(d-1+\sqrt{(d-1)^2-4\mu}\right)=\frac{1}{4} \left(d-1 + |d-1-2\lambda|\right)=\frac{\lambda}{2}\,,
$$
hence the `good' Frobenius solutions for $w(y)$ correspond via \eqref{w} to the `good' Frobenius solutions for $v(y)$. This implies by \eqref{mu} that
for each $\mu$-eigenvalue there is a corresponding eigenvalue
\begin{equation}\label{mu-lambda}
\lambda=\frac{1}{2} (d-1+\sqrt{(d-1)^2-4\mu}).
\end{equation}
This correspondence can be used to get a lower bound on the number of positive eigenvalues. To this end, consider the function
\begin{equation}\label{w0}
w^{(n)}_0(y)=(1-y^2)^{1/2}\, y^{\frac{d-1}{2}}\, v_0^{(n)}(y)=(1-y^2)^{1/2} y^{\frac{d+1}{2}} f'_n(y)\,.
\end{equation}
 This function solves Eq.\eqref{An} for $\mu=d-2$, hence by the Sturm oscillation theorem the number of $\mu$-eigenvalues below $(d-2)$ is equal to the number of zeros of $f_n'(y)$  which, from
\eqref{N}, is $n$ (for $d=3,4$) or $n-1$ (for $d=5,6$). Using  the correspondence \eqref{mu-lambda} between $\mu$-eigenvalues and the eigenvalues $\lambda$, we infer that  $f_n(y)$ has exactly $n$ (for $d=3,4$) or $n-1$ (for $d=5,6$) eigenvalues $\lambda>d-2$.
In addition, for $d=5$ and $n\neq 0$ the function $v_1^{(n)}(y)=(1-y^2)^{-1} y f_n'(y)$, which by \eqref{w} corresponds to $w^{(n)}_0(y)$ for $\lambda=d-2=3$,  is smooth at $y=1$ as follows from \eqref{case1}, hence $\lambda_1^{(n)}=3$ is the eigenvalue for each $n\neq 0$ in $d=5$.

The numerical computations of eigenvalues indicate that, apart from the eigenvalues $\lambda\geq d-2$ discussed above,  in dimensions $d=3,4,5$
there are no additional eigenvalues with positive real part, while for $d=6$ there is exactly  one additional positive real eigenvalue $\lambda_1^{(n)}$. Moreover, for each $f_n$ there are infinitely many negative real eigenvalues (note that the eigenvalues with $\Re(\lambda)<\frac{d-1}{2}$ a priori need not be real). It appears there are no other eigenvalues, confirming   the spectrum \eqref{spectrum}. In the most relevant case here, $n=1$, the numerically computed eigenvalues are displayed in Table~1.

\begin{table}[h]
  \centering
  \begin{tabular}{|>{\centering}m{1ex}| >{\centering}m{9ex} >{\centering}m{3ex} >{\centering}m{9ex} >{\centering}m{9ex} >{\centering}m{9ex} >{\centering\arraybackslash}m{9ex}|}
    \toprule
    $d$ & $\lambda_{1}^{(1)}$ & $\lambda_{0}^{(1)}$ & $\lambda_{-1}^{(1)}$ & $\lambda_{-2}^{(1)}$ & $\lambda_{-3}^{(1)}$ & $\lambda_{-4}^{(1)}$ \\
    \midrule
      3 & 6.33363 & 1 & $-0.51861$ & $-1.75203$ & $-2.88873$ & $-3.97644$ \\
      4 & 3.99883 & 1 & $-0.39021$ & $-1.58542$ & $-2.71468$ & $-3.81626$ \\
      5 & 3 & 1 & $-0.28177$ & $-1.44755$ & $-2.57372$ & $-3.68316$ \\
      6 & 2.42624 & 1 & $-0.17996$ & $-1.30848$ & $-2.41983$ & $-3.52385$ \\
    \bottomrule
  \end{tabular}
  \vskip 0.2cm
  \caption{\small{The six largest eigenvalues of linear modes about $f_{1}$.}}
  \label{tab:eigen}
\end{table}
\vspace{-0.5cm}
 An important consequence of the above considerations is the existence of a unique self-similar solution $f_1(y)$ with exactly one unstable mode.
 This solution is an expected candidate for the critical solution which leads us to:
  \begin{conj*}
   For $d\in\{3,4,5,6\}$ the threshold of blowup is given by the codimension-one stable manifold of the self-similar solution $f_1$.
  \end{conj*}
  In the next section we present
  the numerical corroboration of this conjucture.

\section{Threshold for blowup}
Before discussing  evolution, we show in Figure~1 and Table~2 the numerically computed profiles and parameters of the self-similar solution $f_1$ in  $3\leq d \leq 6$.
\begin{figure}[h]
  \centering
  \includegraphics[width=0.7\textwidth]{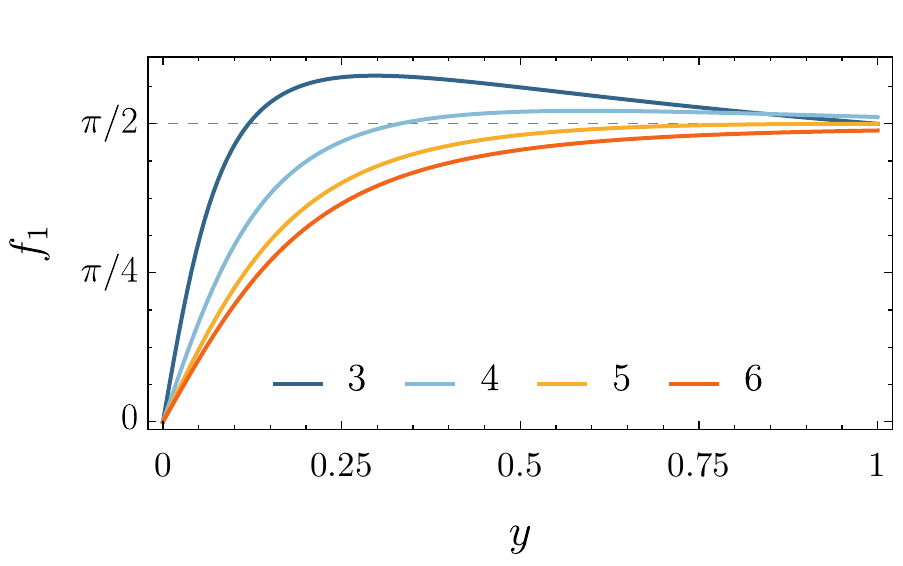}
  \caption{\small{Profiles of  $f_{1}$ in
    $d=3,4,5,6$.}}
  \label{fig:SelfSimilarProfilesF1}
\end{figure}

\begin{table}[ht]
  \centering
  \begin{tabular}{|>{\centering}m{1ex}| >{\centering}m{9ex} >{\centering}m{9ex} >{\centering\arraybackslash}m{9ex}|}
    \toprule
    $d$ & $f_{1}'(0)$ & $f_{1}(1)$ & $f_{1}'(1)$ \\
    \midrule
    3 & 21.75741 & $\pi/2$ & $-0.30566$ \\
    4 & 10.9953 & 1.60634 & $-0.10654$ \\
    5 & 7.82119 & $\pi/2$ & $0$ \\
    6 & 6.71508 & 1.53534 & $0.059052$ \\
    \bottomrule
  \end{tabular}
  \vskip 0.2cm
  \caption{\small{Shooting parameters of  $f_{1}$ in
      $d=3,4,5,6$.}}
  \label{tab:1}
\end{table}

%
Numerical simulations of blowup require special methods that are able to resolve vanishing spatio-temporal scales as the solution develops a singularity. In our previous study of generic blowup for equation \eqref{eqv} we used a moving mesh method \cite{bb}. This method is computationally costly which is a serious drawback in the present context because the precise fine-tuning to the threshold requires many runs of the code.
For this reason we propose here a different method which is based on specially designed similarity-like coordinates
for which self-similar solutions are asymptotic stationary states. We believe that our method is interesting on its own and could be useful in numerical simulations of self-similar blowup for other evolution equations.

Note that standard similarity coordinates $(s,y)$ are not suitable for numerical evolution because the time of blowup is not known a priori which leads to the the gauge mode instability, as described above. To go around this difficulty we shall use a self-correcting coordinate system that adapts
to the upcoming blow-up time as the solution evolves. This new coordinate system $(\tau,\rho)$ is defined as follows
\begin{equation}
  \label{tau_rho}
    r=e^{-\tau}\rho,\quad t=t(\tau),\quad \frac{dt}{d\tau}=e^{-\tau} h(\tau),
\end{equation}
where the function $h(\tau)$, defining the relation between the numerical slow-time $\tau$ and $t$, will be chosen below (for $h=1$ the coordinates $(\tau,\rho)$ coincide with the similarity coordinates $(s,y)$). The new dependent variables $(V,P)$ are defined as follows
\begin{equation}
u(t,r)=V(\tau,\rho),\quad \partial _{t}u(t,r)=e^\tau P(\tau,\rho)\,.
\end{equation}
In terms of these new variables the wave map equation takes the form
\begin{subequations}
  \label{eqnum}
  \begin{align}
    \label{eqV}
    \partial_{\tau} V&=h P-\rho \partial_{\rho} V,\\
    \partial_{\tau} P&=h\left(\partial_{\rho\rho} V+\frac{d-1}{\rho}\partial_{\rho} V-\frac{d-1}{2\rho^2}\sin(2V)\right)-P-\rho \partial_{\rho} P.
  \end{align}
\end{subequations}
Differentiating \eqref{eqV} with respect to $\rho$ and evaluating at $\rho=0$ we get
\begin{equation}\label{eqV0}
  \partial_\tau(\partial_\rho V(\tau,0))+\partial_\rho V(\tau,0)=h(\tau)\partial_\rho P(\tau,0).
\end{equation}
Choosing
\begin{equation}\label{h}
h(\tau)= \left(\partial_{\rho} P(\tau,0)\right)^{-1}
\end{equation}
and solving \eqref{eqV0} we obtain
\begin{equation}\label{Vinf}
  \partial_\rho V(\tau,0)=1+c e^{-\tau}.
\end{equation}
Thus, regardless of whether the solution  blows up or not, its gradient $ \partial_\rho V(\tau,0)$ tends asymptotically to $1$.
The gradient $\partial_{\rho} P(0,\tau)$ also remains bounded but its asymptotic value depends on an endstate of evolution. To see this, note that
\begin{equation}\label{u_tr}
    \partial_{rt} u(t,0)=e^{2\tau}\partial_\rho P(\tau,0)\,,
  \end{equation}
which implies that in the case of dispersion we have $\partial_{\rho} P(0,\tau)\rightarrow 0$ (hence $h(\tau) \rightarrow \infty$), while in the case of blowup along the self-similar solution $f_n$ we have $\partial_{\rho} P(0,\tau)\rightarrow 1/f_n'(0)$ (hence $h(\tau) \rightarrow f_n'(0)$).

The boundedness of the gradient of any solution is a very desirable
feature of our formulation since it allows us to solve the system
(25), with $h$ given by \eqref{h}, using standard finite difference
methods on a uniform grid. More specifically, we use a fourth order centered finite difference scheme to
approximate spatial derivatives in the interior of the numerical grid
$0<\rho<\rho_{\textrm{max}}$. Near the origin $\rho=0$ we use
symmetries of functions $V$ and $P$ to evaluate derivative stencils,
while near the artificial boundary $\rho_{\textrm{max}}$ we use
one-sided schemes. We evolve this semi-discrete system in time with
fifth order adaptive Runge-Kutta method, known as DOPRI5
\cite{Hairer}. To suppress spurious high frequencies we add
standard dissipation terms. As an outer boundary
of the radial grid we typically take $\rho_{\text{max}} \approx 2 f_n'(0)$,
where $f_n\left(\frac{r}{T-t}\right)$ is the expected self-similar
endstate of evolution. This guarantees that the grid includes the past
light cone of the singularity.
\begin{rem*}
  Our numerical method can be viewed as a simplified moving mesh
  method combined with a Sundman transformation.  In a moving mesh
  method each mesh point can move independently from other points (the
  motion of mesh points is governed by a prescribed mesh density
  function).  In our case, the mesh points form a more rigid structure
  as reflected by a simple relation between $r$ and $\rho$.  The
  relative scale of $r$ and $\rho$ is governed only by a single degree
  of freedom $h(\tau)$.  The same parameter also dictates the relative
  scales of the time variables $t$ and $\tau$ just as a Sundman
  transformation would.
\end{rem*}
We illustrate our numerical results for initial data of the form
\begin{equation}\label{id}
V(0,\rho)=\frac{A \rho}{\cosh{\rho}} = P(0,\rho)\,.
\end{equation}
\addtolength{\topmargin}{-1.2pc}
\addtolength{\textheight}{2.4pc}
For large $A$ we find that $h(\tau)\rightarrow f_0'(0)$ (which corresponds to generic  blowup governed by the self-similar solution $f_0$), while for small $A$ we find that $h(\tau)\rightarrow \infty$ (which corresponds to dispersion to zero). Using bisection we fine tune the amplitude $A$ to the critical value $A_{*}$ with precision of 32 digits. For such marginally critical amplitudes we observe that for intermediate times $h(\tau)$ approaches  $f_1'(0)$. This is illustrated in Fig.~2 which shows the marginally sub- and supercritical evolutions in $d=6$ (the plots for $d=3,4,5$ look very similar).
\begin{figure}[h]
  \centering
  \includegraphics[width=0.85\linewidth]{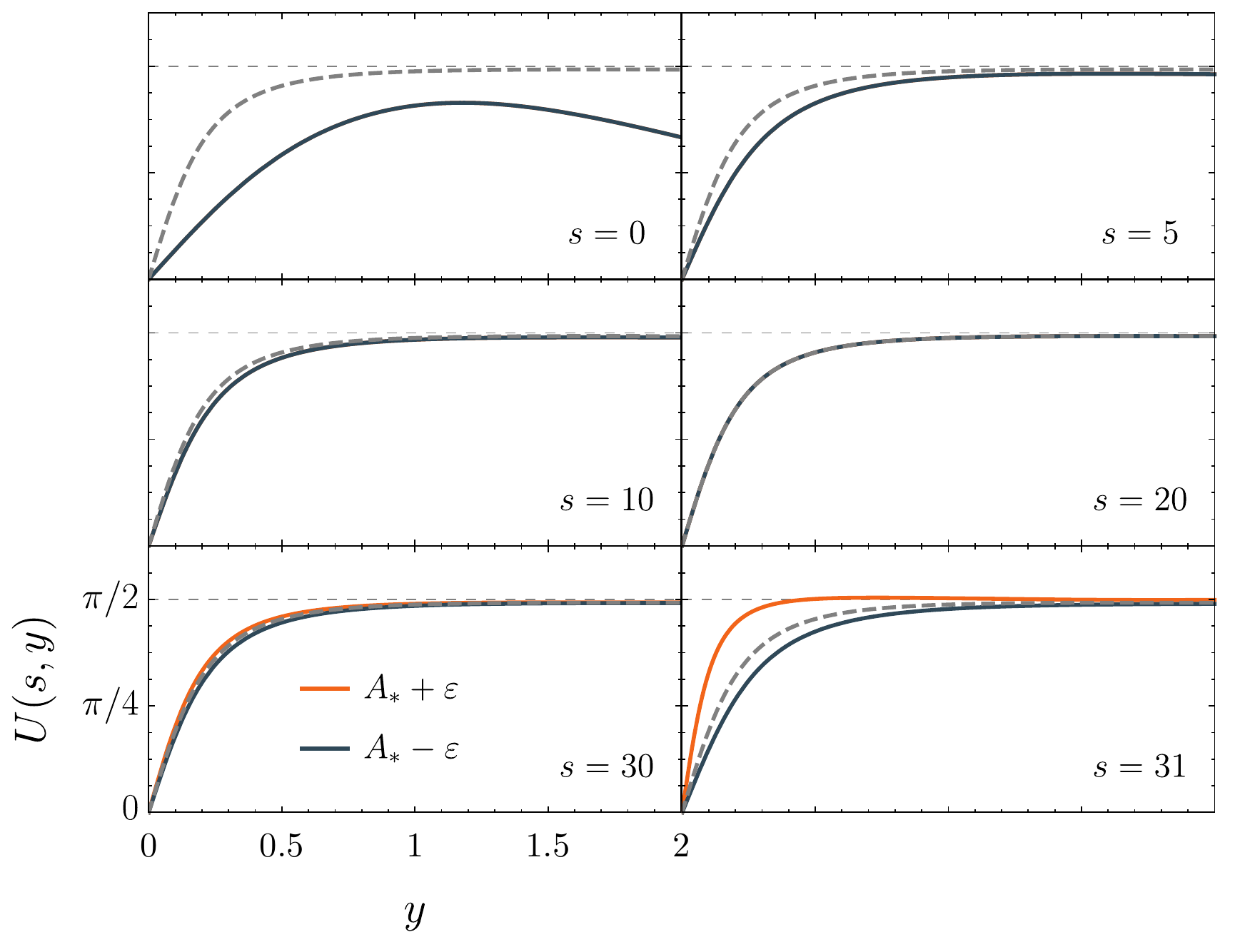}
  \caption{\small{The evolution of a pair of initial data \eqref{id} in $d=6$ with marginally sub- (blue line) and supercritical (red line) amplitudes $A=A_*\pm \varepsilon$, where $\varepsilon$ is of the order $10^{-32}$. The solutions evolve together, approach the intermediate attractor $f_1$ (dashed line), and eventually depart from it in opposite directions.}}
  \label{fig:snapshots_d6}
\end{figure}

\noindent To compare the results with the predictions of the linear perturbation analysis from section~3, we now
  translate the results into the similarity coordinates $(s,y)$. To this end we need to determine the blowup time $T$.  This is done as follows. First, we integrate equation  $\frac{dt}{d\tau}=e^{-\tau} h(\tau)$ to get $t(\tau)$. For intermediate times $t(\tau)$ develops a plateau which yields a rough estimate for $T$. Having that, we compute
$e^{-s}=T-t$ and then
\begin{equation}\label{UtoV}
  \partial_{y} U(s,0)\approx e^{\tau-s}\partial_{\rho} V(\tau,0).
\end{equation}
From the linear perturbation analysis it follows that for intermediate times (when the solution is close to the threshold) the left hand side of \eqref{UtoV} is well approximated by
\begin{equation}\label{interm}
  \partial_{y} U(s,0)\approx f_{1}'(0)+a_{1}e^{\lambda_{1}^{(1)}s}+a_{0}e^{s}+a_{-1}e^{\lambda_{-1}^{(1)}s},
\end{equation}
\nopagebreak
where the coefficient $a_1 \sim A-A_*$ is very small.
\newpage

Since our estimate of the blowup time is not precise, this approximation involves the gauge mode instability with a nonzero coefficient $a_0(T)$.
 Fitting the formula \eqref{interm} to  the right hand side of \eqref{UtoV}, we get the coefficients  $a_{1}$, $a_{0}$, and $a_{-1}$, which depend on the estimated value of $T$. Finally, performing bisection with respect to $T$ we determine the precise blowup time $T_{*}$ for which $a_{0}(T_{*})=0$.

 The result of such a fit for the marginally critical evolution from Figure~2 is shown in Figure~3 (to plot both the
sub- and supercritical solutions against the same variable
$s=-\log(T-t)$, the blow-up time $T$ was chosen to be the average of
$T_{*}$ for the sub- and supercritical solutions). The fit shows excellent agreement with the results of the linear perturbation analysis which makes us feel confident that our conjecture is true.

\begin{figure}[h]
  \centering
  \includegraphics[width=0.9\linewidth]{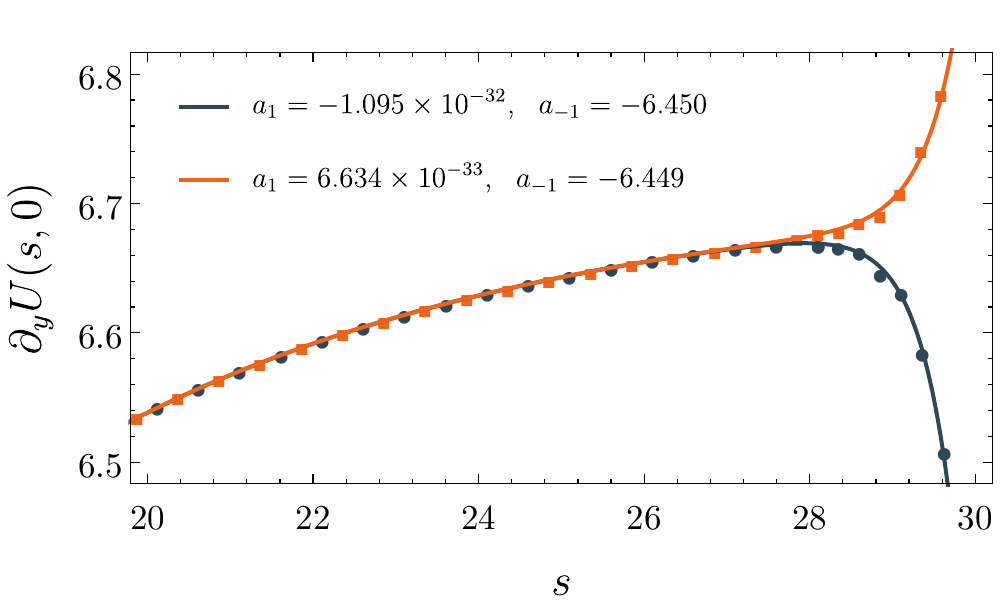}
  \caption{\small{We plot the gradients at the origin for the marginally critical solutions shown in
    Figure~2.  The gradients converge to $f_{1}'(0)\approx6.71$ with a rate given
    by the stable eigenvalue $\lambda_{-1}^{(1)}\approx -0.18$ and
    almost the same coefficients $a_{-1}$.  At a later time, the gradients
    separate and stray away from $f_{1}'(0)$, each growing with the
    rate given by the unstable eigenvalue
    $\lambda_{1}^{(1)}\approx2.43$ and very small coefficients
    $a_{1}$ of opposite signs.}}
  \label{fig:fitd6}
\end{figure}

\emph{Acknowledgement.}  This research was supported in part by the Polish National Science Centre grant no. DEC-2012/06/A/ST2/00397. We gratefully acknowledge the support of the Alexander von Humboldt Foundation.

\end{document}